\documentclass[12pt]{amsart}
\newtheorem{thm}{Theorem}[section]

\newtheorem{defn}[thm]{Definition}
\newtheorem{lem}[thm]{Lemma}
\newtheorem{cor}[thm]{Corollary}
\newtheorem{rem}[thm]{Remark}

\usepackage{amsmath}
\usepackage{amsxtra}
\newcommand{\z}{\mathcal Z}
\newcommand{\gr}{{\rm gr}\ }
\newcommand{\id}{{\rm id}}
\usepackage{amscd}
\usepackage{amsthm}
\usepackage{amsfonts}
\usepackage{amssymb}
\usepackage{eucal}

\newcommand{\g}{{\mathfrak{g}}}
\newcommand{\gl}{{\mathfrak{gl}}}
\newcommand{\n}{{\mathfrak{n}}}
\renewcommand{\sl}{{\mathfrak{sl}}}
\newcommand{\ch}{{\rm ch}}
\newcommand{\C}{{\mathbb C}}
\newcommand{\F}{{\mathcal F}}
\newcommand{\Z}{{\mathbb Z}}

\begin{document}
\title[Fusion products and Kostka polynomials] {Fusion products,
  cohomology of $GL_N$ flag manifolds and Kostka polynomials}
\date{Submitted: 9/2003. Accepted for publication in IMRN, 12/2003}
\author{Rinat Kedem}
\begin{abstract}
  This paper explains the relation between the fusion product of
  symmetric power $\sl_n$ evaluation modules, as defined by Feigin and
  Loktev, and the graded coordinate ring $R_\mu$ which describes the
  cohomology ring of the flag variety $\mathcal F l_{\mu'}$ of $GL_N$.
  The graded multiplicity spaces appearing in the decomposition of the
  fusion product into irreducible $\sl_n$-modules are identified with
  the multiplicity spaces of the Specht modules in $R_\mu$.  This
  proves that the Kostka polynomial gives the character of the fusion
  product in this case.  In the case of the product of
  fundamental evaluation modules, we give the precise correspondence
  with the reduced wedge product, and thus the usual wedge space
  construction of irreducible level-1 $\widehat{\sl}_n$-modules in the
  limit $N\to \infty$. The multiplicity spaces are $W(\sl_n)$-algebra
  modules in this limit.
\end{abstract}

\maketitle
\section{Introduction}
The graded tensor product (fusion product) of finite-dimensional
evaluation representations of the Lie algebra $\g\otimes \C[t]$, where
$\g$ is a simple Lie algebra, was introduced by Feigin and Loktev in
\cite{FL}. The usual tensor product of irreducible finite-dimensional
$\g$-modules $\{\pi_{\mu^{(i)}}, i=1,...,N\}$, with highest weights
$\mu^{(i)}$, decomposes into the direct sum of irreducible
$\g$-modules:
\begin{equation}\label{decom}
\pi_{\mu^{(1)}}\otimes\cdots \otimes \pi_{\mu^{(m)}} \simeq
\underset{\lambda}{\oplus} K_{\lambda,\boldsymbol\mu}\ \pi_\lambda 
\end{equation} 
($\boldsymbol \mu \overset{\rm def}{=} (\mu^{(1)},...,\mu^{(m)})$),
where the integers $K_{\lambda,\boldsymbol\mu}$, 
products of Clebsch Gordan
coefficients, are the multiplicities of $\pi_\lambda$ in the tensor
product. 

The fusion product of \cite{FL} is a $\g$-equivariant
grading on this product, so that the multiplicities
$K_{\lambda,\boldsymbol\mu}$ become the
generating functions for the dimensions of the graded components of
the multiplicity space.  It was conjectured \cite{FL} that these
polynomials are related to the Kostka polynomials in the special cases
where such polynomials are defined.

On the other hand, in the case where $\mu^{(i)}=\mu_i \omega_1$ are
the highest weights of the symmetric power representation, and
$\mu=(\mu_1,...,\mu_m)$ with $\sum_i \mu_i=N$, it is known
\cite{GP,CP} that Kostka polynomials give the graded dimensions of the
cohomology ring of the partial flag variety $\mathcal F l_{\mu'}$ of
$GL_N$. The Weyl group, which is the symmetric group $S_N$, acts on
this ring, preserving degree.  As a graded space and as an
$S_N$-module, the cohomology ring is isomorphic to the coordinate ring
$R_\mu$, which is a certain finite-dimensional quotient of
$\C[z_1,...,z_N]$ (see section 3 for the precise definition) on which
$S_N$ acts by permutation of variables. The decomposition of $R_\mu$
into irreducible $S_N$-modules is
$$
R_\mu \simeq \underset{\lambda\vdash N}{\oplus} W_\lambda \otimes
M_{\lambda,\mu}
$$
where $W_\lambda$ are the Specht modules, and $M_{\lambda,\mu}$ is
a graded multiplicity space.  The graded multiplicity of
$M_{\lambda,\mu}$ is a Kostka polynomial.

In this paper, we use a simple argument to show how the Schur-Weyl
duality allows us to identify $M_{\lambda,\mu}$ with the multiplicity
spaces of irreducible $\sl_n$-modules in the graded tensor product of
\cite{FL}, in the special case where all evaluation modules are
symmetric power representations.

In \cite{FJKLM}, a much more general definition of a level-restricted,
generalized Kostka polynomial was made in terms of the characters of
the level restricted fusion product, which is a quotient of the fusion
product of \cite{FL} by a certain ideal.  In the limit where level
$k\to \infty$ this gives a definition of the unrestricted, generalized
Kostka polynomial.  In the case where $\mu^{(i)}$ are highest weights
of the form $\mu_i \omega_j$ with $\omega_j$ fundamental weights, it
is still a conjecture that this definition overlaps with that of
\cite{Sh,SW}.

The method presented here for the symmetric power representations does
not rely on explicit formulas for the Kostka polynomial, but it relies
on the calculation of \cite{GP} for the character of the quotient ring
$R_\mu$ described in Section 3.

In the special case where $\pi_{\mu^{(i)}}$ are all fundamental
representations, the ring $R_\mu$ appears also in the more usual
construction of the (reduced) wedge product of evaluation
representations. These products appear in the well-known construction
of the basic (level-1, irreducible) $\widehat{\sl}_n$-modules in terms
of the inductive limit of the wedge product, modulo the action of the
Heisenberg algebra. The reduced wedge product is the finite
dimensional version of this, the wedge product modulo the action of
the Heisenberg algebra. It is also possible to construct integrable
$\widehat{\sl_n}$-modules by taking the appropriate inductive limit of 
special graded tensor products \cite{FF}.  We explain the relation
between the reduced wedge product and the graded tensor product of the
fundamental representations.

The paper is organized as follows. In section 2, we recall the
definition of the graded tensor product of \cite{FL}. In section 3,
the construction of $R_\mu$ due to Garsia and Procesi \cite{GP} is
explained. Section 4 gives a construction of the graded tensor product
of the fundamental representation of $\sl_n$ in terms of the action of
the symmetric group and the Schur-Weyl duality.
In section 5, we explain the analogous
construction for the reduced wedge product. Section 6 explains the
limit $N\to \infty$ of this special fusion product.  In section 7, we
explain the case of the product of symmetric power representations of
$\sl_n$.

The author is pleased to acknowledge useful discussions with
N. Reshetikhin and E. Stern.

\section{The graded tensor product}
The graded tensor product was introduced by Feigin and Loktev in
\cite{FL} for any simple Lie algebra $\g$. It is called the ``fusion
product'' in that reference because it is motivated by the fusion
product of affine algebra modules in conformal field theory, in the
case where the level is generic.

Let $\widehat{\g}\simeq \g\otimes \C[t^{-1},t]]\oplus \C c$
 be the affine algebra associated with $\g$, with central element
$c$.  On an irreducible
$\widehat{\g}$-module, $c$ acts as a constant $k\in \C$ (the level).
If $k$ is generic, the Weyl module $W_\mu = {\rm
  Ind}_\g^{\g\otimes t^{-1}\C[t^{-1}]} \pi_\mu$, induced from the
finite dimensional irreducible $\g$-module $\pi_\mu$ by the action of
$\g\otimes t^{-1} \C[t^{-1}]$, (where $\g\otimes t\C[t]$ acts
trivially on $\pi_\mu$ and $c$ acts as $k$) is an irreducible
level-$k$ module. In this case, $\pi_\mu$ is called the {\em top component}
of $W_\mu$.

Let $\{z_1,...,z_N\}$ be distinct finite points in $\C P^1$.  By a
$\widehat{\g}$-module at the point $z_i$, we mean the module
on which the element $x\otimes f(t)$ acts as $x\otimes f(t-z_i)$.
There exists a level-$k$ action of $\widehat{\g}$ on the tensor
product of $N$ level-$k$ modules $W_{\mu^{(i)}}$ at the
points $z_i$. The loop algebra $\g\otimes\C[[t,t^{-1}]$ acts by the
usual coproduct with the action on $W_{\mu^{(i)}}$ being at
$z_i$, while cocycle on the product is the sum of the canonical
cocyles.  This action has the same level $k$ as the modules
$W_{\mu^{(i)}}$, and is referred to as the fusion action. The product
of $\widehat{\g}$-modules with this action is called the fusion
product. We denote it by $W_{\boldsymbol \mu}(\mathcal Z)$ ($\mathcal
Z = (z_1,...,z_N)$).

Let $\mathcal A$ denote space of regular functions on $\C
P^1\setminus\{z_1,...,z_N\}$ vanishing at infinity, with possible
poles at $z_i$. Then $\g\otimes
\mathcal A$ acts on $W_{\boldsymbol \mu}(\mathcal Z)$ (with trivial
central charge) and the space of $\g\otimes \mathcal A$-invariant
functionals on $W_{\boldsymbol \mu}(\mathcal Z)$ is isomorphic to the
dual of the top component, $\pi_{\mu^{(1)}}\otimes \cdots\otimes
\pi_{\mu^{(N)}}$ in $W_\mu(\z)$, for generic $k$. Its dual space is
the space of coinvariants.

On the other hand, the fusion action of $\g\otimes \C[t]$ (in fact, of
$\n_-[t]$) on the tensor product of highest weight vectors
$v_{\mu^{(i)}}\in\pi_{\mu^{(i)}}$ generates the top component
$\pi_{\mu^{(1)}}\otimes \cdots\otimes \pi_{\mu^{(N)}}$ if $z_i$ are
distinct. Moreover, there is a natural grading on
$U(g\otimes\C[t])v_{\mu^{(1)}}\otimes\cdots\otimes v_{\mu^{(N)}}$, by
degree in $t$ (see below for a precise description of the grading),
which is $\g$-equivariant. Therefore, this gives a graded version of
the tensor product of the finite dimensional representations, which is
described by the formula (\ref{decom}) with
$K_{\lambda,\boldsymbol{\mu}}$ graded multiplicities.  This is the
idea behind the definition of the fusion product introduced in
\cite{FL}.  Let us recall the precise definition.

Given $\mu=(\mu_1,\mu_2,...,\mu_n)$, a partition with $\mu_1\geq
\cdots\geq \mu_n\geq 0$, $\mu_i\in\Z$, denote by $\pi_{\mu}$ the irreducible
$\sl_n$-module with $\sl_n$-highest weight
$(\mu_1-\mu_2)\omega_1+\cdots+(\mu_{n-1}-\mu_n)\omega_{n-1}$ where
$\omega_i$ are the fundamental weights.

Let $z$ be a formal parameter. 
Denote by $\pi_{\mu}[z] \simeq \pi_{\mu}\otimes \C[z]$ the
{\em formal evaluation module} of
$\sl_n[t]=\sl_n\otimes\C[t]$, which acts in the usual way:
$$
(x\otimes f(t))\circ (w\otimes g(z)) = x\ w \otimes f(z)g(z),
 \qquad w\in\pi_\mu, \quad x\in\sl_n.
 $$
 All tensor products are over $\C$.  The formal evaluation module
 is infinite-dimensional. The evaluation module at a complex number
 $a$ is the image of the evaluation map ${\rm id}\otimes \phi_a:
 \pi_\mu\otimes \C[z]\to \pi_\mu$ corresponding to evaluation of the
 formal parameter $z$ at the complex number $a$. The evaluation module
 at a complex number is finite-dimensional. It is denoted in the same
 way, $\pi_\mu[a]$, and it will be made clear in what follows which type of
 evaluation module is under consideration.

Given a sequence of $N$ evaluation representations $\{\pi_{\mu^{(i)}}
 [z_i],\
i=1,...,N \}$ 
corresponding to partitions $\mu^{(i)}$, let
 $\boldsymbol\mu=(\mu^{(1)},...,\mu^{(N)})$ and define
$ V_{\boldsymbol\mu}({\mathcal Z})$ to be the tensor
product:
\begin{eqnarray}
V_{\boldsymbol\mu}(\z)&=&
\pi_{\mu^{(1)}}[z_1]\otimes\cdots\otimes \pi_{\mu^{(N)}}[z_N]
\nonumber \\
&\simeq&
\left(\pi_{\mu^{(1)}}\otimes\cdots\otimes
\pi_{\mu^{(N)}}\right)\otimes \C[z_1,...,z_N].
\end{eqnarray}
It contains a subspace 
$$\F_{\boldsymbol \mu}(\mathcal Z)= U\left(\n_-\otimes \C[t]\right)\ 
\left((v_{\mu^{(1)}}\otimes \cdots \otimes v_{\mu^{(N)}})\ \otimes
  1\right),$$
where $v_{\mu^{(i)}}$ is the highest weight vector of
$\pi_{\mu^{(i)}}$ and $1\in\C[z_1,...,z_N]$. The action is by the
usual co-product on evaluation modules. When $z_i$ are taken to be
formal parameters, the space $\F_{\boldsymbol \mu}(\z)$ is graded by
homogeneous degree in the $z_i$'s. If, instead, $z_i$ are specialized
to be distinct
complex numbers, we can define a filtration of the space as follows.
Let
$$
\F_{\boldsymbol\mu}^{(m)}(\z) = U^{(m)}(\n_-\otimes \C[t])
v_1\otimes v_2 \otimes \cdots \otimes v_N,
$$
where $U^{(m)}$ is spanned by homogeneous elements of
degree $m$ in $t$. The subspaces
$$
\F_{\boldsymbol\mu}^{(\leq n)}(\z)=\oplus_{m\leq n}
\F_{\boldsymbol\mu}^{(m)} (\z)
$$
define a filtration of $\F_{\boldsymbol\mu}(\z)$.

The fusion product can be defined in either setting.
In the case where $z_i$ are specialized to distinct complex numbers,
the graded tensor product of \cite{FL} is the following:
\begin{defn} The graded tensor product of evaluation modules at the distinct
  complex numbers $\{z_1,...,z_N\}$ is
$$
\F_{\boldsymbol\mu}^* \overset{\rm def}{=} \gr \F_{\boldsymbol\mu}(\z).
$$
\end{defn}

It is possible to define the graded tensor product in the setting
where $z_i$ are formal variables also. In this setting, the space
$\F_\mu(\z)$ is an infinite-dimensional graded space. Define the subspace
$$
\widetilde{\F}_{\boldsymbol\mu}^{(n)} (\z)= (\C[z_1,...,z_N]\otimes_\C
{\F_{\boldsymbol\mu}}^{(\leq n-1)})(\z) \cap
{\F_{\boldsymbol\mu}}^{(n)}(\z)\subset \F_\mu^{(n)}(\z).
$$
Then the graded tensor product can be defined as
\begin{defn}\label{fusion} The graded tensor product of formal
  evaluation modules is
$$
\F_{\boldsymbol\mu}^{*} = 
\underset{n
  \geq 0}{\oplus}
\mathcal{F}_{\boldsymbol\mu}^{(n)}(\z)/
\widetilde{\F}_{\boldsymbol\mu}^{(n)}(\z). 
$$
\end{defn}

The two definitions are clearly equivalent. In either case, the space
$\F_{\boldsymbol\mu}^{*}$
is a finite-dimensional $\sl_n$-module, of dimension equal
to the dimension of the tensor product of the finite dimensional
$\sl_n$-modules $\pi_{\mu^{(i)}}$. Since the action of $\sl_n$
does not change the degree, the graded components are
$\sl_n$-modules. The decomposition of the graded tensor product into
irreducible $\sl_n$-modules,
$$
\F_{\boldsymbol\mu}^{*} \simeq \underset{\lambda}{\oplus}
\pi_\lambda \otimes \mathcal M_{ \lambda,{\boldsymbol\mu}},
$$
defines the multiplicity spaces $\mathcal M_{ \lambda,{\boldsymbol\mu}}$.
The multiplicity spaces are graded spaces, and the purpose of this
paper is to give an explicit description of these spaces.

It was conjectured in \cite{FL} that the graded multiplicities,
$$
\ch_q \mathcal M_{ \lambda,{\boldsymbol\mu}} = \sum_i q^i \dim
\mathcal M_{ \lambda,{\boldsymbol\mu}}[i]
$$
(for any graded space $V$ we denote by $V[i]$ its graded component)
are related to the generalized Kostka polynomials of \cite{Sh,SW}.  In
the paper \cite{SW}, the generalized Kostka polynomial
was defined for the
case where $\mu^{(i)}$ are partitions of the form $(j)^m$
corresponding to rectangular Young diagrams, or highest weights $m
\omega_j$.  One can expect to relate the $q$-dimension of $\mathcal M_{
  \lambda,{\boldsymbol\mu}}$ to the generalized Kostka polynomial in
this case, although this is still a conjecture.

Assuming that this conjecture is correct, one can take this character
to be the definition of a more general Kostka polynomial. This was the
subject of \cite{FJKLM}, where the level-restricted version of the
most general Kostka polynomial was defined as a quotient of the graded
tensor product. Here, we discuss only the generic level case, which 
is obtained from \cite{FJKLM} by taking the limit as the level $k$
tends to infinity.

\begin{rem} The level restriction is related to the fact that the
  decomposition of the fusion product of irreducible representations
  of $\widehat{\g}$, in the case where $k$ is an integer, is
  determined by the Verlinde rule and not by the usual Clebsch-Gordan
  rule. In that case the dimensions $K_{\lambda, \boldsymbol \mu}$ are
  smaller. We do not consider this case in this paper.
\end{rem}

\section{Polynomial $S_N$-representations and the Kostka polynomials}\label{GP}
Here, we review some useful facts from \cite{GP} about the
decomposition of certain quotients of the space of polynomials in $N$
variables into irreducible symmetric group modules. These quotients
were obtained \cite{GP,CP} as a description of the cohomology
ring of the flag variety of $GL_N$. The symmetric group is the Weyl
group of $GL_N$ so it acts on the cohomology ring. As a symmetric
group module, the ring corresponding to the flag $\mathcal F l_{\mu'}$
is isomorphic to $\C S_N/S_\mu$ where $S_\mu$ is the Young subgroup.

Let $z_i$ be formal variables, and consider the action of $S_N$ on the
ring $\C[z_1,...,z_N]$ by permutation of variables. Let $\mu$ be a
partition of $N$, $\mu = (\mu_1,...,\mu_m)$ with $\mu_1\geq \mu_2 \geq
\cdots \geq \mu_m>0$. Let $\{a_1,...,a_m\}$ be distinct complex numbers,
and let $\mathcal X = (x_1,x_2,...,x_N)$ denote the $N$-tuple of complex
numbers with the first $\mu_1$ variables equal to $a_1$, the next
$\mu_2$ equal to $a_2$, etc. 

Denote by $I(\mathcal X)$ the ideal of functions in $\C[z_1,...,z_N]$
which vanish under the evaluation map 
\begin{eqnarray}
\phi_{\mathcal X}&: \C[z_1,...,z_N]&\to \C \nonumber\\
& f(z_1,...,z_N)&\mapsto f(x_1,...,x_N)\label{evaluationmap}.
\end{eqnarray}
This ring is invariant under the action of
Young subgroup of $S_N$ given by the partition $\mu$, $S_\mu =
S_{\mu_1}\times S_{\mu_2} \times \cdots \times S_{\mu_m}$. Here, by
$S_{\mu_1}$ we mean the permutation group of the first $\mu_1$
variables $(z_1,...,z_{\mu_1})$ and so forth.

Let $I_\mu$ be the intersection of the ideals in the $S_N$-orbit of
$I(\mathcal X)$:
$$
I_\mu = \cap_{\sigma\in S_N} I(\sigma{\mathcal X}),
$$
where $S_N$ acts by permutation of the indices in
$(x_1,...,x_N)$. Functions in $I_\mu$ vanish at all points in the
$S_N$-orbit of $(x_1,...,x_N)$.

The ideal $I_\mu$ is $S_N$-invariant by definition, hence the quotient
ring 
\begin{equation}\label{amu}
A_\mu = \C[z_1,...,z_N]/I_\mu
\end{equation}
is an $S_N$-module, isomorphic to the representation of $S_N$ acting
on left cosets of $S_\mu$. The space $A_\mu$ is filtered by the
homogeneous degree in the variables $z_i$. Define $R_\mu$ to be the
associated graded space:
$$
R_\mu = \gr A_\mu
$$
 We give an explicit combinatorial description
of $R_\mu$ below for completeness, but we will use it only in the
special case when $\mu = (1)^N$, which is explained first.

In the case $\mu=(1)^N$, let $J_N$ be the ideal in the ring
$\C[z_1,...,z_N]$, generated by symmetric polynomials in $N$
variables, of degree greater than 0:
 For example, the generators of
this ideal can be taken to be the elementary symmetric functions,
$E_1(z_1,...,z_N),...,E_N(z_1,...,z_N)$, where
$$
E_m(z_1,...,z_N) = \sum_{1\leq i_1<i_2<\cdots<i_m\leq N} z_{i_1}\cdots z_{i_m}.
$$
\begin{equation}\label{JN}
J_N = \langle E_1(z_1,...,z_N),...,E_N(z_1,...,z_N)\rangle.
\end{equation}
Since $J_N$ is invariant with respect to
the action of $S_N$, $S_N$ acts on the quotient space
\begin{equation}\label{Rn}
R_N:= R_{(1)^N}=\C[z_1,...,z_N]/J_N.  
\end{equation}
This is a finite-dimensional vector space. In fact,
\begin{lem}\cite{GP}
  As a symmetric group module, $R_N$ is is isomorphic to the regular
  representation:
$$
R_N \simeq \oplus_{\lambda\vdash N} W_{\lambda}\otimes M_{\lambda, (1)^N}
$$
where $W_{\lambda}$ is the irreducible $S_N$ Specht (left-) module,
corresponding to the partition $\lambda$, and $M_{\lambda,(1)^N}$ is a
multiplicity space, of dimension equal to the dimension of
$W_\lambda$.
\end{lem}

Note that the ideal $J_N$ is generated by homogeneous polynomials, so
the quotient $R_N$ is a graded space.  The symmetric group action
preserves homogeneous degree in $z_i$, so the Specht modules are
spanned by polynomials of fixed degree.  Thus, the multiplicity spaces
$M_{\lambda,(1)^N}$ are graded spaces
and the $q$-character (or $q$-dimension)
$$
\ch_q(M_{\lambda,(1)^N}) = \sum_i q^i \dim M_{\lambda,(1)^N}[i]
$$
is the
Kostka polynomial $\widetilde{K}_{\lambda, (1)^N}(q)$, where 
\begin{equation}\label{kostilde}
\widetilde{K}_{\lambda,\mu}(q)=q^{n(\mu)}K_{\lambda,\mu}(1/q).
\end{equation}
Here, $n(\mu)=\sum_i (i-1)\mu_i$.  The normalization factor ensures
that this is a polynomial in $q$, as it should be: the lowest degree
polynomial in $R_N$ is 1, of degree 0.  Here, $K_{\lambda,\mu}(q)$ is the
usual Kostka polynomial \cite{Mac}, defined, for example, as the
transition matrix between Schur polynomials and Hall-Littlewood
polynomials:
$$
S_\lambda(\mathbf x) = \sum_{\mu\leq \lambda} K_{\lambda,\mu}(q)
P_\mu(\mathbf x,q).
$$

For completeness, we recall here the description of the coordinate ring
$R_\mu$ \cite{GP} in the general case.  The ring $R_\mu$ is the associated
graded space of the quotient $A_\mu$. For a function
$f(z_1,...,z_N)\in I_\mu$, define the top term to be the term of
highest homogeneous degree in $\{z_i\}$. Such terms generate an ideal
$J_\mu$, and are obtained by taking the limit of $f(z_1,...,z_N)$ as
$x_i\to 0$ for all $i$. The ring $R_\mu$ is the quotient of
$\C[z_1,...,z_N]$ by the ideal $J_\mu$. It has an explicit
combinatorial description as follows.

Let $d_k(\mu)= N-(\mu'_1 + \cdots + \mu_{N-k}')$, where $\mu'$ is the
conjugate partition to $\mu$.  Define the set
$$C_\mu(z_1,...,z_k) = \{E_r(z_1,...,z_k): k-d_{k}<r\leq k\}.$$
If
$d_{k}(\mu)=0$ the set is empty. If $k=N$, the set includes all
elementary symmetric polynomials of positive degree in $N$ variables.
Let $C_\mu$ be the set of polynomials in the $S_N$-orbit of
$C_\mu(z_1,...,z_k)$ for all $k$, 
$$C_\mu=\underset{\underset{k>0}{\sigma\in S_N}}
{\cup}\sigma C_\mu(z_1,...,z_k).$$
For example, the set
$C_{(1^N)}$ is just the case $k=N$ above, and the set $C_{(N)}$
includes all polynomials other than $1$.

\begin{thm}\cite{GP}
  The ideal $J_\mu$ is generated by $C_\mu$. The graded space $R_\mu$
  is the quotient ring
$$
 R_\mu = \C[z_1,...,z_N]/J_\mu .
$$
\end{thm}
The
decomposition of $R_\mu$ into irreducible $S_N$-modules,
$$
R_\mu \simeq \oplus_{\lambda\vdash N} W_\lambda \otimes M_{\lambda,\mu},
$$
gives the graded multiplicity spaces $M_{\lambda,\mu}$.
\begin{thm}\cite{GP}
$$
{\rm ch}_q\ M_{\lambda,\mu} = \widetilde{K}_{\lambda,\mu}(q).
$$
\end{thm}

$R_\mu$, as a representation of the symmetric group, is isomorphic to
the representation of the symmetric group acting on left cosets of the
Young subgroup $S_\mu$, where $S_\mu$ is the stabilizer of $\mathcal X$.

\section{Graded tensor product fundamental $\sl_n$-modules}\label{cc}
The simplest example of the graded tensor product product of
$\sl_n$-representations is the $N$-fold product of the fundamental
$n$-dimensional evaluation representation associated with
$\pi=\pi_{(1)}$ of highest weight $\omega_1$.  We have $\mu^{(i)} =
(1)$ for all $i\in \{1,...,N\}$, and
$$
V_N(\z) :=V_{(1^N)}(\mathcal Z) = \pi[z_1]\otimes\cdots\otimes
\pi[z_N] \simeq
\pi^{\otimes N}\otimes_\C \C[z_1,...,z_N].
$$
We consider here the tensor product of formal evaluation modules.

This case is interesting for two reasons.  First, it is isomorphic to
the reduced finite-dimensional wedge product, which is recalled in the
next section. Second, the inductive limit as $N\to\infty$ gives a
construction of the basic representations of $\widehat{\sl}_N$
\cite{FL}.  In terms of the cohomology ring of the flag variety, its
character is related to the complete flag.

Denote by $v_0,...,v_{n-1}$ the standard basis of $\pi\simeq \C_n$,
where $v_0$ is the highest weight vector.  The subspace $\F_N:=
U(\n_-\otimes\C[t]) w_0$, with $w_0=(v_0\otimes\cdots\otimes
v_0)\otimes 1\in \pi^{\otimes N}\otimes \C[z_1,...,z_N]$, can be
characterized largely in terms of the action of the symmetric group
$S_N$. This is the approach taken in this section.

As a symmetric group module, the space $V_N$ is a tensor product of
two $S_N$-modules, say, $V_N=\Pi\otimes G$, where $\Pi= \pi^{\otimes
  N}$, on which $S_N$ acts by permutation of factors in the tensor
product, and $G=\C[z_1,...,z_N]$, on which $S_N$ acts by permutation
of variables as before. 
We consider the {\em diagonal left action} of symmetric group on
this product. That is, for any $\sigma\in S_N$,
\begin{eqnarray*}
&& \sigma \circ \left((v_{i_1}\otimes v_{i_2}\otimes \cdots)\otimes
f(z_1,z_2,...)\right) = \\
& & (v_{\sigma(i_1)}\otimes v_{\sigma(i_2)}\otimes \cdots) \otimes
f(\sigma(z_1),\sigma(z_2),\cdots)). 
\end{eqnarray*}

By definition, $\sl_n\otimes \C[t]$ acts on the tensor product by the
coproduct, which commutes with the diagonal left action of
$S_N$. In addition, the cyclic vector $w_0$ is $S_N$-invariant.
Therefore, we have that
$$\F\subset ( \pi^{\otimes N}\otimes\C[z_1,...,z_N])^{S_N}.$$
(This is
true for any $N$-fold product of identical representations, not only
the fundamental one.) Here, if $W$ is a vector space on which a group $G$
acts, by $W^{G}$ we mean the subspace on which $G$ acts trivially.
Note that the $S_N$-action preserves homogeneous degree.

Next, consider the definition \ref{fusion} of the graded tensor
product.  We can choose an $S_N$-invariant basis for the
representatives of the quotient $\F^{(n)}_N/
\widetilde{\F}^{(n)}_N$.  To see this, let $p^{(n)}\in
\F^{(n)}_N$ and let $\overline{p}^{(n)}$ be a representative of $p^{(n)}$
in the quotient. That is, there exists an
expression
\begin{equation}\label{residues}
p^{(n)} = \overline{p}^{(n)} + \sum_i\sum_{m<n}
f_{i}^{(n-m)}(z_1,...,z_N) p^{(m)}_{i}
\end{equation}
where $p_i^{(m)}\in \F^{(m)}_N$, and 
$f_{i}^{(n-m)}(z_1,...,z_N)$ are homogeneous
polynomials of degree $n-m$.
Since both $p^{(n)}$ and each of the vectors $p_i^{(m)}$ are
$S_N$-invariant, we have for any $\sigma\in S_N$:
$$\overline{p}^{(n)} -\sigma(\overline{p}^{(n)}) = \sum_i\sum_{m<n}
(f_{i}^{(n-m)}(\mathbf z)-f_{i}^{(n-m)}(\sigma(\mathbf z))) p_i^{(m)}
\in \widetilde{\F}_N^n.$$
Thus, the difference is zero modulo $\widetilde{\F}_N^n$, so it vanishes in
the quotient. Therefore, we can take as a basis for the
coset representatives $S_N$-invariant vectors $\overline{p}^{(n)}$.

This means that in choosing an $S_N$-invariant basis,
the functions $f_{i}^{(n-m)}$ appearing in (\ref{residues}) are
symmetric functions. That is,
$$ \widetilde{\F}^{(n)}_N \simeq
(\C[z_1,...,z_N]^{S_N} \otimes \mathcal{F}_N^{(\leq n-1)})\cap
\mathcal{F}_N^{(n)}.
$$
Since the graded space $\F_N^{(n)}$ has degree $n$, such
polynomials must be of positive degree. 

This means that $\widetilde{\F}^{(n)}_N$ is in contained in the ideal
in $\F_N$ generated by symmetric polynomials of positive degree.
Defining the ideal $\mathcal J$ by $\F^{*}_{(1)^N}= V_N^{S_N}/\mathcal
J$, we have shown that $\mathcal J\subset (\pi^{\otimes N}\otimes
J_N)^{S_N}$, where $J_N$ is the ideal defined in (\ref{JN}).

We will show the equality: The graded component of degree $n$ of the
ideal in $\F_N$ generated by symmetric functions of positive degree
contains $\widetilde{\F}_N^{(n)}$, $\mathcal J\cap \F_N^{(n)}\subset
\widetilde{\F}^{(n)}_N.$  This follows from Lemma
\ref{vspan}, which shows that any vector in $V_N(\z)^{S_N}$ can be
expressed as a linear combination (with coefficients symmetric
functions) of vectors in $\F_N$.  That is, $V_N(\z)^{S_N}\simeq
\F_N\otimes_\C \C[z_1,...,z_N]^{S_N}$.  It follows that if $p^{(n)}\in
\F_N^{(n)}$ and $p^{(n)}=\sum_{i=0}^{n-1} f_{n-i} v_{i}$ with
$f_{n-i}$ symmetric functions of degree $n-i>0$, then $v_i\in
\F_N^{(\leq i)}\otimes \C[z_1,...,z_N]^{S_N}$ and hence
$p^{(n)}\in\widetilde{\F}_N^{(\leq n-1)}$.

\begin{lem}\label{vspan}
The space $V_N(\z)^{S_N}$ is generated by $\F_N$ as a
$\C[z_1,...,z_N]^{S_N}$-module. 
\end{lem}
\begin{proof}
  First, note that in the case of the fundamental module, $\pi$ is
  generated by the action of the commuting elements 
  $$\{b_i= E_{1,i+1}, i=1,...,n-1\}$$
  acting on $v_0$, with
  $b_j v_0 = v_j$.  Moreover, $b_i b_j v_0 = 0$
  for any $i,j$.
  
  We denote the vector $v_{i_1}\otimes\cdots\otimes v_{i_N}\in
  \pi^{\otimes N}$ by
  $[i_1,...,i_N]$.  The space $V_N(\z)^{S_N}$ is clearly spanned by elements
  of the form
\begin{equation}\label{v}
v_{\mathbf N}^{(\mathbf m)} =
\sum_{\sigma\in S_N} \sigma
 [\underbrace{n-1,...,n-1}_{N_{n-1}},\underbrace{n-2,...,n-2}_{N_{n-2}},..., 
\underbrace{0,...,0}_{N_0}]\otimes z_{\sigma(1)}^{m_1}\cdots
  z_{\sigma(N)}^{m_N}.
\end{equation}
Here, $\mathbf m$ is an $N$-tuple of integers, and $\mathbf N$ is the
$n$-tuple $(N_{n-1},...,N_0)$.

If $m_j=0$ for all $j>N-N_0$, then $v_{\mathbf N}^{(\mathbf
  m)}\in \F_N$.
This is because $\F_N$ is spanned by elements of the form
$$
(b_{i_1}\otimes t^{m_1}) \cdots (b_{i_k}\otimes t^{m_k}) w_0  =
\sum_{\sigma\in S_N} \sigma[i_1,...,i_k,0,...,0]\otimes
z_{\sigma(1)}^{m_1}\cdots z_{\sigma(k)}^{m_k}
$$
with $n-1\geq i_1\geq i_2\geq \cdots \geq i_k\geq 1$.
Thus in general, let $\nu(\mathbf m,\mathbf N)$ denote the number of non-zero
components $m_j$ with $j>N-N_0$.

The lemma is proven by induction on $\nu(\mathbf m,\mathbf N)$. If
$\nu(\mathbf m,\mathbf N)=0$ then $v_{\mathbf N}^{(\mathbf m)}\in
\F_N$.  Suppose $\nu(\mathbf m,\mathbf N)>0$, we will show that
$v_{\mathbf N}^{(\mathbf m)}$ can be expressed as a sum of the form
\begin{equation}\label{vvv}
\sum_{\mathbf m'} f_{\mathbf m'}(z_1,...,z_N) v_{\mathbf N}^{(\mathbf m')}
\end{equation}
where $f_{\mathbf m'}(z_1,...,z_N)$ are symmetric functions, and $\nu(\mathbf
m', \mathbf N))<\nu(\mathbf m,\mathbf N)$ for all $\mathbf m'$
occuring in the sum.  Given this
fact, it follows by induction that we can rewrite any element in
$V_N(\z)^{S_N}$ in the form (\ref{vvv}) with $\nu(\mathbf m',\mathbf
N)=0$, proving the lemma.

Define a partition of the set of integers $\{1,...,N\}$ into subsets $J_{0},
...,J_{n-1}$, with $J_{n-1} = \{1,...,N_{n-1}\}, J_{n-2} =
\{N_{n-1}+1,...,N_{n-1}+N_{n-2}\}$ and so forth.
Denote the corresponding Young subgroup
of $S_N$ by $S_{\mathbf J} = S_{J_0}\times S_{J_1}\times\cdots\times
S_{J_{n-1}}$. By $\mathbf m(J_i)$, denote the $N_{i}$-tuple composed
of the elements $m_j$ where $j\in J_i$.

Let $M(J;\mathbf m(J))$ be the monomial symmetric function in the
variables $\{z_i: i\in J\}$:
$$
M(J;\mathbf m(J)) = \sum_{\sigma\in S_J} \prod_{i\in J} z_{\sigma(i)}^{m_i}.
$$
Up to a constant multiple, the expression (\ref{v}) for $v_{\mathbf
  N}^{(\mathbf m)}$ is equal to
\begin{equation}
\sum_{\sigma \in S_N/S_{\mathbf J}}
\sigma [n-1,...,0]\otimes
 \prod_{i=0}^{n-1} M(\sigma(J_i);\mathbf m(J_i)),
\label{rewritea} 
\end{equation}
because $\sigma\in S_{\mathbf J}$ acts trivially on the vector
$[n-1,...,0]$ and on $\prod M(J_i,\mathbf m(J_i))$.
 The summation is over the left coset representatives.

One can always rewrite
$$M(J_0;\mathbf m(J_0)) = F_1-F_2$$
where
$$
F_1 = ((N_1+\cdots +N_{n-1})!)^{-1}\sum_{\sigma\in S_N}
\prod_{i\in J_0} z_{\sigma(i)}^{m_i}
$$
is a completely symmetric function. Since $M(J_0;\mathbf m(J_0))$
is symmetric in each set of variables $\{z_j: j\in J_i\}$ separately,
so is $F_2$. But each monomial in $F_2$ has at least one non-zero
power of some $z_i$ with $i\in J_k$ with $k\in \{1,...,N-1\}$. This is
because all of the terms containing only $z_j$ ($j\in J_0$) are
all contained in $F_1$.  Therefore,
\begin{equation}\label{ftwo}
F_2 = \sum_{\mathbf m'} c_{\mathbf m'}\prod_{i=0}^{n-1} M(J_i, \mathbf m'(J_i))
\end{equation}
with some constants $c_{\mathbf m'}$, where $\mathbf m'(J_0)$
contains fewer non-zero elements than $\mathbf m(J_0)$.

Using this fact to rewrite (\ref{rewritea}), we have
the difference
of the two kinds of terms:
\begin{eqnarray*}
& & \sum_{\sigma\in S_N/S_J} \sigma[n-1,...,0]\otimes
\sigma\left(\prod_{i=1}^{n-1} M(J_i,\mathbf m(J_i))(F_1-F_2)\right) \\
&=&
F_1 \sum_{\sigma\in S_N/S_J} \sigma[n-1,...,0]\otimes
\sigma\prod_{i=1}^{n-1} M(J_i,\mathbf m(J_i))\\
& & \quad - 
\sum_{\sigma\in S_N/S_J} \sigma[n-1,...,0]\otimes
\sigma \prod_{i=1}^{n-1} M(J_i,\mathbf m(J_i)) F_2.
\end{eqnarray*}
Here, since $F_1$ is completely symmetric, it can be factored out of
the summation, since $S_N$ acts on it trivially. The first term
is proportional to
\begin{equation}\label{firstterm}
\sum_{\tau\in S_N}\prod_{j\in J_0}z_{\tau(j)}^{m_j}
\sum_{\sigma\in S_N/S_{\mathbf J}} \sigma[
n-1,...,0]\otimes \prod_{i=1}^{n-1} M(\sigma(J_i);\mathbf m(J_i)),
\end{equation}
Let $m'_j=m_j$ if $j\in J_1\cup\cdots\cup J_{n-1}$, and $m'_j=0$ if
$j\in J_0$. Then 
(\ref{firstterm}) is equal to
$$
f(z_1,...,z_N) v_{\mathbf N}^{\mathbf m'},
$$
where $f$ is a symmetric function. This term is in the
$\C[z_1,...,z_N]^{S_N}$-span of $\mathcal F_N(\z)$.

The second term involves the function
$$
\prod_{i=1}^{n-1} M(J_i,\mathbf m(J_i)) F_2,
$$
which is again symmetric in each set of variables $\{m_j, j\in
J_i\}$ for each $i$. Due to (\ref{ftwo}), the second term is therefore
expressible as a sum of several terms of the same form (\ref{ftwo})
with $\mathbf m''$ replacing $\mathbf m'$, and where the number of
non-zero elements in $\mathbf m''(J_0)$ is the same as in $\mathbf
m'(J_0)$. It is strictly less than $\nu(\mathbf m,\mathbf N)$.

The lemma follows by induction.
\end{proof}

\begin{rem}
  The lemma simply states that $V_N(\z)^{S_N} = U(\gl_n\otimes \C[t])
  w_0$.
\end{rem}

The lemma means that any element of degree $i$ in $V_N(\z)^{S_N}$ which is
zero modulo symmetric functions of positive degree is also zero modulo
$\F_N^{(<i)}$. Therefore, $\F^{*}_{(1)^N}\simeq V_N(\z)^{S_N}/(\pi^{\otimes
  N}\otimes J_N)^{S_N}$ as $\sl_n$-modules and as graded vector
spaces.  As a consequence, 
\begin{thm}
$$\F^{*}_{(1)^N} \simeq (\pi^{\otimes N}\otimes R_N)^{S_N}.$$
\end{thm}
Recall that as $S_N$-modules
$$
R_N \simeq {\underset{\lambda\vdash N}{\oplus}}W_\lambda\otimes
M_{\lambda, (1)^N},
$$
where $M_{\lambda,\mu}$ is a graded multiplicity space and
$W_\lambda$ is the irreducible Specht module.  Also, by the Schur-Weyl
duality, the action of $S_N$ on $\pi^{\otimes N}$ centralizes the
action of $\sl_n$, and the tensor product decomposes as
$$
\pi^{\otimes N}\simeq 
\underset{{\mu\vdash N, \ell(\mu)\leq n}}{\oplus} \pi_\mu \otimes W_\mu,
$$
where $\pi_\mu$ is the irreducible representation of $\sl_n$, and
$W_\mu$ is an irreducible $S_N$-module.

Putting all these facts together,
$$
\left( (\underset{\mu}{\oplus}\pi_\mu \otimes W_\mu)\otimes
  (\underset{\lambda}{\oplus} W_\lambda\otimes
  M_{\lambda,(1)^N})\right)^{S_N} \simeq \underset{\mu,\lambda}{\oplus}
\pi_\mu\otimes M_{\lambda,(1)^N} \otimes \left( W_\mu\otimes
  W_\lambda\right)^{S_N}.
$$
The coefficient of the trivial representation in the tensor product
$W_\mu\otimes W_\lambda$ is equal to 1 if and only if $\lambda = \mu$.
Thus, we have
\begin{thm}\label{together}
  As graded vector spaces and $\sl_n$-modules, the graded tensor
  product of the basic representations is equal to
$$
\F^{*}_{(1)^N} \simeq \underset{\underset{\ell(\mu)\leq n}{\mu\vdash
    N}}{\oplus} \pi_\mu\otimes M_{\mu,(1)^N},
$$
where $M_{\mu,(1)^N}$ is the graded multiplicity space of the
irreducible $S_N$-representation $W_\mu$ in the quotient ring
$R_N$.
\end{thm}
Finally in terms of characters,
\begin{cor}\label{ccchar}
$$
ch_q \F^{*}_{(1)^N} = \sum_{\mu\vdash N, \ell(\mu)\leq n} ch(\pi_\mu)
\widetilde{K}_{\mu,(1)^N}(q). 
$$
\end{cor}

\section{The reduced wedge space}
The $N$-fold graded tensor product of the fundamental representations
is closely related to the wedge product of evaluation representations
with formal evaluation parameters. Recall the wedge product
construction of the basic representations of the affine Lie algebra
$\widehat{\sl}_n$: They are obtained as the quotient of the
semi-infinite wedge product, stabilized at infinity, by the image of
a Heisenberg algebra, which can be described as the central
extension of the loop algebra of the center of $\gl_n$, acting on the
wedge product.

This suggests a finite-dimensional analogue of the $N$-fold wedge
product, which is the quotient of the finite-dimensional wedge product
by the image of the same Heisenberg algebra, acting with central
charge zero.  We call this the reduced wedge product.  In
the appropriate limit $N\to\infty$, this reduced wedge product becomes
the irreducible level-1 $\widehat{\sl}_n$-module.\footnote{The construction
explained here appeared previously in a preprint by the author and E. Stern
(unpublished).}

We show here that the reduced wedge product and the graded tensor
product of Section \ref{cc} are isomorphic, as $\sl_n$-modules and
as graded vector spaces.

Recall that there are $n$ irreducible, integrable representations of
$\widehat{\sl}_n$ with level $k=1$, $L(\Lambda_i)$, with highest
weight vector of weight $\Lambda_i$, the fundamental affine weight,
with $i=0,...,n-1$. Let us recall the wedge space construction of
$L(\Lambda_i)$.

Let $\pi$ be the fundamental $\sl_n$-module as before, considered now
as a $\gl_n$-module, with basis
$v_0,...,v_{n-1}$ with $v_0$ the highest weight vector. Let $z$ be a
formal variable, and $\pi(z)\simeq \pi\otimes \C[z,z^{-1}]$ the
evaluation module of $\widehat{\gl}_n$.  The space
$\pi(z)\wedge\pi(z)$ can be realized in $\pi(z)\otimes\pi(z)$
as the span of vectors of the form
\begin{eqnarray*}
v_i z^m \wedge v_j z^n &=& v_i z^m \otimes v_j z^n -
v_j z^n\otimes v_i z^m \\ &=& v_i\otimes v_j z_1^m z_2^n - v_j\otimes v_i
z_1^n z_2^m,\quad 0\leq i,j\leq n-1; \ n,m\in \Z,
\end{eqnarray*}
where $z_i$ refers to the parameter $z$ appearing in the $i$-th factor
in the tensor product.  Therefore if we again consider the diagonal
action of the symmetric group on the space $\pi^{\otimes
  N}\otimes\C[z_1^{\pm1},...,z_N^{\pm1}]$, the $N$-fold wedge product
of fundamental evaluation representations can be realized as the
subspace on which the symmetric group $S_N$ acts as the alternating
representation. The algebra $\widehat{\gl}_n$ acts on the $N$-fold
wedge product by the usual co-product, with central charge $0$.

Denote the basis of $\pi(z)$ by $v_{i}^k=v_i\otimes z^k$. 
The semi-infinite
wedge space $F^{(j+1,k)}$ ($n-1\geq j\geq 0, k\in \mathbf Z$) is the $\C$-span of all vectors which differ from
\begin{equation}\label{span}
w_{j+1}^k=v_{j}^k\wedge v_{j-1}^k\wedge\cdots \wedge v_{0}^k\wedge
v_{n-1}^{k+1}\wedge v_{n-2}^{k+1}\wedge\cdots\wedge v_{0}^{k+1}\wedge
v_{n-1}^{k+2}\wedge \cdots
\end{equation}
in a finite number of factors. For sufficiently large $N$,
a basis vector of $F^{(j+1,k)}$ is different from (\ref{span}) only in the
first $N$ factors.

There is a level-1 action of the algebra $\widehat{sl}_n$ and the
Heisenberg algebra $\mathcal H$ generated by $\{B_n = {\rm id} \otimes
t^n,\ n\neq 0\}\subset \widehat{\gl}_n$ on this space, and they
centralize each other.  For any $k$ and $j$, the vector $w_j^k$ is a
highest weight vector with respect to both algebras.  The irreducible
$\widehat{sl}_n$ representation $L(\Lambda_{\overline{i}})$ is
isomorphic to the quotient of $F^{(i,k)}$ (with any integer $k$ and
$\overline{i}=i\ {\rm mod \ n}$) by the image of the Heisenberg
algebra $\mathcal H_- =\{{\rm id} \otimes t^{-n}\ n>0\}$ acting on
$F^{(i,k)}$:
$$
L(\Lambda_{\bar i}) \simeq F^{(i,k)}/U(\mathcal H_-) F^{(i,k)}.
$$

The stabilization requirement for $F^{(i,k)}$ can be realized as
follows. Let $F^{(i,k)}_m$ be the subspace of $F^{(i,k)}$ spanned by
vectors which may differ from (\ref{span}) only in the first $n m + i$
factors.  The space is spanned by the elements
$$
v_{i_1}^{k_1}\wedge v_{i_2}^{k_2}\wedge\cdots \wedge
v_{i_{nm+i}}^{k_{nm+i}} \wedge v_{n-1}^{k+m+1}\wedge
v_{n-2}^{k+m+1}\wedge \cdots
$$
where $0\leq i_j\leq n-1$ and $k_j<k+m+1$.
We have
$$
F^{(i,k)}_{0}\subset F^{(i,k)}_{1} \subset\cdots\subset F^{(i,k)}_m
\subset \cdots
$$
where the inductive limit gives the full semi-infinite wedge space. 

It follows that $F^{(i,k)}_m = (\wedge^{(N)} \pi(z)) \wedge
w_{n-1}^{k+1+m}$, where $N=mn+i$. (Strictly speaking, this wedge
product makes sense when considering the full Clifford module
$F=\oplus F^{(i,k)}$ and the action of the Clifford algebra on it.)
Thus, the space $F^{(i,k)}_m $ is isomorphic to the space generated by
the action of $\gl_n\otimes \C[t^{-1}]\simeq
\sl_n\otimes\C[t^{-1}]\oplus\mathcal H_-$ on a cyclic vector of the
finite $N$-fold wedge product obtained by truncating the highest
weight vector of the semi-infinite
wedge product $F^{(i,k)}$ after the $N^{\rm th}$ factor. For
example, in the case $k=-m$, this vector is
\begin{equation}\label{hwv}
\overline{w}_i^{-m} 
=
v_{i-1}^{-m}\wedge v_{i-1}^{-m}\wedge\cdots v_0^{-m} \wedge v_{n-1}^{-m+1}
\wedge \cdots \wedge v_0^{-m+1} \wedge \cdots \wedge v_{n-1}^{0}
\wedge \cdots
\wedge  v_{0}^{0},
\end{equation}
and $U(\gl_n\otimes\C[t^{-1}]) \overline{w}_i^{-m}
 = \wedge^{mn+i} \pi[z^{-1}]$.
The isomorphism with $F^{(i,-m)}_m$ is the wedge product with the
highest weight vector of $F^{(n,1)}$, $w_{n-1}^{1}$. Note that the
choice of $k$ is immaterial as the difference is only in overall
normalization of the character (or equivalently, multiplication by an
overall monomial in $z_i^{\pm1}$) and we choose a
$k$-independent convention for the normalization of the character
as follows.

To be consistent with the grading used in Section \ref{cc}, we define
the $q$-character of $F^{(i,k)}_m$ by taking the degree of $t$ to be
$1$. This means
$$
\ch_q F^{(i,k)}_m = \sum_{j\leq 0} q^j \dim(F^{(i,k)}_m[j])
$$
where we normalize the highest weight vector for any $k,i$ to have
degree 0, and $t$ (or $z_i$) has degree 1. This is a polynomial in
$q^{-1}$, which we compute below.

The irreducible module $L(\Lambda_i)$ is the quotient of the
semi-infinite wedge product with respect to the image of $\mathcal H_-$.
The vector $w_{n-1}^{k+1+m}$
is not in the image of $\mathcal H_-$, so the
quotient of $F^{(i,k)}_m$ by the image of the Heisenberg algebra is
obtained by taking the quotient of the finite wedge product by the
image of $\mathcal H_-$, then taking the wedge product with
$w_{n-1}^{k+1+m}$. That is,
$$
F_m^{(i,k)}/U(\mathcal H_-)F_m^{(i,k)} \simeq \wedge^N\pi[z^{-1}]/U(
  \mathcal H_-)(\wedge^N\pi[z^{-1}])
$$
(where the isomorphism is multiplication by a certain monomial in $z_i$,
and wedging with $w_{n-1}^{k+1+m}$), and
$$
\left(\lim_{m\to\infty}F_m^{(i,k)}\right)/U( \mathcal H_- ) F^{(i,k)}=
\lim_{m\to\infty}\left( F_m^{(i,k)}/U( \mathcal H_- )F_m^{(i,k)} \right).
$$
  Therefore, one can reverse the order of taking the
limit $m\to\infty$ and taking the quotient by the image of $\mathcal H_-$.

On the vector $\overline{w}_i^{-m}$ (and on any vector in the finite
wedge product), $U(\mathcal H_-)$ acts as multiplication by symmetric
polynomials in $N$ variables $z_i^{-1}$ of positive degree.  Thus, the
quotient of the finite wedge product with cyclic vector (\ref{hwv})
can be described as follows. Let
$\pi[z^{-1}]=\pi\otimes \C[z^{-1}]$.
The $N$-fold wedge product of such modules
is an infinite-dimensional space (recall $z_i$ are formal variables)
\begin{equation}\label{wedge}
\wedge^N \pi[z^{-1}] \simeq (\pi^{\otimes N}\otimes \C[z_1^{-1},...,z_N^{-1}])
^A
\end{equation}
where the superscript $A$ means the subspace on which the diagonal
action of $S_N$ is as on the alternating representation: $\sigma\in
S_N$ acts by $(-1)^{{\rm sgn}(\sigma)}$.

The reduced wedge space is a quotient of this space by the image of
the Heisenberg algebra $\mathcal H_-$. But the generator $B_a=
\id\otimes t^{-a}$ ($a>0$) acts on the wedge product by multiplication
by $z_1^{-a} + \cdots + z_N^{-a}$. Such polynomials form a basis for
the space of symmetric functions in $N$ variables of positive degree,
so they generate the ideal $J_N(\z^{-1})$ (i.e. with $z_i$ replaced by
$z_i^{-1}$) in $\C[z_1^{-1},...,z_N^{-1}]$.  The image of the
Heisenberg algebra in (\ref{wedge}) is the ideal generated by
multiplication by elements in $J_N(\z^{-1})$.  As in the previous
section, we have
$$
\wedge^N \pi[z^{-1}] /{\rm Im}\ {\mathcal H_-} \simeq (\pi^{\otimes N}\otimes
R_{N}(\z^{-1}))^A,
$$
where $R_N(\z^{-1})$ is the ring (\ref{Rn}) with $z_i$ replaced by $z_i^{-1}$.
Again, we use the theorem of \cite{GP} and the Schur-Weyl duality to
re-write this as
$$
\wedge^N \pi[z^{-1}] /{\rm Im}\ {\mathcal H_-} \simeq
\underset{\underset{\lambda\vdash N}{\mu\vdash N, \ell(\mu)\leq n}}{\oplus}
 \pi_\mu \otimes M_{\lambda, (1)^N}(\z^{-1}) \otimes (W_\mu\otimes
W_\lambda)^A .
$$
The summand is non-zero if and only if $\mu=\lambda'$
(the conjugate partition of $\mu$). Thus we have
\begin{lem}
$$
\wedge^N \pi[z^{-1}] /{\rm Im}\ {\mathcal H_-}\simeq \underset{\mu:
  \ell(\mu)\leq n}{\oplus}
\pi_\mu\otimes M_{\mu',(1)^N}(\z^{-1}).
$$
\end{lem}
Therefore,
\begin{eqnarray*}
\ch_q\ \wedge^N \pi[z^{-1}] /{\rm Im}\ {\mathcal H_-} &=& \sum_{\mu\vdash N,
\ell(\mu)\leq n} \ch(\pi_\mu) q^{-N(N-1)/2} K_{\mu',(1)^N}(q)\\ & =&
\sum_{\mu\vdash N,
\ell(\mu)\leq n}  \ch(\pi_\mu) K_{\mu,(1)^N} (q^{-1}),
\end{eqnarray*}
where we used a well-known identity for the Kostka polynomial. 

\begin{lem}
\begin{equation}\label{char}
\ch_q F^{(i,k)}_m/{\rm Im}\ \mathcal H_-  = q^{n(\mu_0')}
 \sum_{\underset{\mu\vdash mn+i}{\ell(\mu)\leq n}} \ch(\pi_\mu) K_{\mu,(1)^N}(q^{-1})
\end{equation}
with $\mu'_0=((n)^m,i)$, $N=mn+i$ and $k\in \Z$.
\end{lem}
\begin{proof}
  Apart from the argument above, the only additional information is
  the normalization, which is determined by the convention that the
  highest weight vector of $F^{(i,k)}$ should have degree 0.  This
  highest weight vector is in the $\sl_n$-module with highest weight
  $w_{i}$, if $i<n$, or the trivial one if $i=n$.  This corresponds to
  $\pi_{\mu_0}$ with $\mu_0' = ((n)^m,i)$ in the decomposition
  (\ref{char}). On the other hand the Kostka polynomial
  $K_{\mu,(1)^N}(q^{-1})$ behaves as a power series in $q^{-1}$ of the
  form $q^{-n(\mu')}(1+O(q^{-1}))$. The lemma follows.

\end{proof}

Thus, up to an overall normalization factor, we see that the
characters of the reduced wedge space and the graded tensor product
(Corollary \ref{ccchar}) are the same. 

The difference in the normalization factor is due to the fact that the
cyclic vector, which is normalized in each case to have degree $0$, is
different. The wedge space is generated by the action of $\sl_n\otimes
\C[t^{-1}]$ on the highest weight vector with respect to the
$\widehat{sl}_n$ action, whereas the graded tensor product is
generated by the action of $\sl_n\otimes \C[t]$ on a cyclic vector
which, in the limit $N\to\infty$, is deep inside the representation --
it is an extremal vector at infinity, in the language of \cite{FS}.

\section{The limit $N\to\infty$}

It was explained in \cite{FL,FF} that the irreducible level-1
$\widehat{\sl}_n$-modules can be obtained as the inductive limit of the
graded tensor product. The same is true for the reduced wedge
product, by construction. Let us consider the behavior of the
characters.  In this section, we change conventions from $q$ to
$q^{-1}$ in the character (so that the character is a polynomial or
power series in $q$ not $q^{-1}$) since this is the standard
convention in the literature.

In the limit $N$ (or $m$) $\to \infty$, for any finite-dimensional
$\pi_\mu$, there is a well-defined limit
$$
\lim_{m\to\infty} q^{-n(\mu_0')} K_{\mu,(1)^{nm+i}}(q).
$$
This limit is a character of the algebra which centralizes the
action of $\sl_n$ on the irreducible level-1 module -- the $W$-algebra
of $\sl_n$ \cite{FKRW}. That is, as an $\sl_n\otimes W(\sl_n)$-module, 
$$
L(\Lambda_i)\simeq \underset{\mu}{\oplus} \pi_\mu \otimes M_\mu,
$$
where $\mu$ is a partition of length less than $n$, such that
$|\mu|\equiv i \mod n$, and $M_\mu$ is
an irreducible $W(\sl_n)$-module whose character is described below.

The $W$-algebra of $\gl_n$ is the central extension of the algebra of
differential operators on $\C[t,t^{-1}]$ (including multiplication by
polynomials in $t$).  These act as symmetric differential operators on
the $N$-fold tensor product of evaluation modules $\pi(z)$, acting
only on $\C[z_1^{\pm1},...,z_N^{\pm1}]$ and not on the vector space.
Since such differential operators commute with the action of the
symmetric group, they also act on the wedge product.  Therefore they
act on the multiplicity space $M_{\mu',(1)^N}$, where the Heisenberg
subalgebra in $W$ acts trivially by definition.  The statement that
the Kostka polynomials tend to the specialized character of an
irreducible $W$-algebra module in the limit
$m\to\infty$ means that in the limit, the $W$ algebra acts irreducibly
on the multiplicity space. This action in the limit has central charge
$n$.

The character formulas of \cite{FKRW} are obtained by using the hook
formula for the Kostka polynomial (Another approach to computing them
is the functional realization, see \cite{FJKLM2} for the case of
$\sl_2$):
\begin{equation}\label{hook}
K_{\mu,(1)^N} = q^{n(\mu')}\frac{(q)_N}{\prod_{x\in\mu}(1-q^{h(x)})},
\end{equation}
where $h(x)$ is the hook length of the box $x$ in the Young diagram of
shape $\mu$: if $x$ has coordinates $i,j$, 
$$h(x)= \mu_i-j + \mu'_j-i+ 1.$$
Here we use the standard notation $(q)_m=\prod_{i=1}^m(1-q^i)$.
It is not difficult to see that
$$
\prod_{x\in\mu}(1-q^{h(x)})= {\prod_{i=1}^n(q)_{\mu_i+i}}
{\prod_{1\leq i<j\leq n}(1-q^{\mu_i-\mu_j+j-i})^{-1}}
$$
if the length of $\mu$ is $n$.

In lemma \ref{char}, $|\mu|$ depends on $m$. However for any module
$\pi_\mu$ which is finite-dimensional in the limit $m\to\infty$, the
parameter $\overline{m}=(\mu_0)_n-\mu_n<\infty$ is finite. In this
case, $\pi_\mu\simeq \pi_{\overline{\mu}}$ as an $\sl_n$-module, where
$\overline{\mu}$ is finite in the limit $m\to\infty$, with parts
$$
\overline{\mu}_j = \mu_j - (m-\overline{m}).
$$
Then
$$
q^{-n(\mu_0')} K_{\mu,(1)^N}(q) =
q^{n(\overline{\mu}')-\overline{m}(i-(\overline{m}-1)n/2)}
(q)_{mn+i}
\frac{\prod_{1\leq i<j\leq n}(1-q^{\overline{\mu}_i-\overline{\mu}_j+j-i})}
{\prod_{i=1}^n(q)_{\mu_i+i}}.
$$
In the limit $m\to\infty$ this becomes:
$$
\lim_{m\to\infty}q^{-n(\mu_0')} K_{\mu,(1)^N}(q) =
\frac{q^{n(\overline{\mu}')-\overline{m}(i-(\overline{m}-1)n/2)}}
{(q)_\infty^{n-1}}\prod_{1\leq i<j\leq
  n}(1-q^{\overline{\mu}_i-\overline{\mu}_j+j-i}). 
$$

In \cite{FKRW}, the irreducible $W(\sl_n)$-modules $M_\mu$ are normalized so that
their specialized charater is
$$
\chi_\mu(q) = q^{n(\mu')+|\mu|}\frac{\prod_{1\leq i<j\leq
    n}(1-q^{\mu_i-\mu_j+j-i})}{(q)_\infty^{n-1}}
$$
(where we factored out the character of the Heisenberg subalgebra
$1/(q)_\infty$). With this convention, we have
$$
\ch_q L(\Lambda_i) = \sum_{\mu: \mu_n=0, |\mu|\equiv i\mod n}
\ch(\pi_\mu) q^{-n(\overline{\mu}_0')} \chi_\mu(q),
$$
where $\overline{\mu}_0' = ((n)^{\overline m},\bar{i})$ depends on $\mu$.

\section{Tensor product of symmetric-power $\sl_n$-modules}

The construction of section \ref{cc} turns out to be sufficient to
deduce the structure of the graded tensor product of any set of
symmetric power evaluation modules, parameterized by highest weights
$\mu_i \omega_1$ ($\mu_i\in \Z_{>0}$). In this section we denote
$\boldsymbol\mu := \mu=(\mu_1,...,\mu_m)$, and we consider the fusion
product of evaluation modules at a complex number.

Recall that any finite-dimensional irreducible $\sl_n$-module
$\pi_\lambda$ is isomorphic to the image of a Young symmetrizer
$c_{t(\lambda)}\in \C S_N$ acting on $\pi^{\otimes|\lambda|}$, where
$\pi$ is the fundamental module.  In general, given a Young diagram
$\lambda$, and choosing any standard tableau $t(\lambda)$ of
shape $\lambda$ on the letters $1,...,|\lambda|=\sum_i \lambda_i$, let
$\mathcal R_{t(\lambda)}\in S_N$ be the stabilizer of the rows of the tableau,
and $\mathcal C_{t(\lambda)}$ be the stabilizer of the columns. The Young
symmetrizer is
$$
c_{t(\lambda)} = \sum_{\sigma\in \mathcal C_{t(\lambda)}}
{\rm sign}(\sigma) \sigma
\sum_{\sigma'\in \mathcal R_{t(\lambda)}} \sigma'. 
$$

Similarly, the evaluation module $\pi_\lambda[a]$ is isomorphic to the
image of the composite map $c_{t(\lambda)}\otimes \phi_a$ acting on
$\otimes_{i=1}^{|\lambda|}\pi[z_i]$, where $\phi_a$ is the evaluation
of the polynomial at the point $z_i=a$ for all $i$ (since the
isomorphism is independent of the particular tableau $t(\lambda)$, we will
drop the tableau notation).  We use this to construct the general
tensor product of symmetric power evaluation modules, and explain how
the graded tensor product can be obtained as a quotient of this.

Let $V_N(\mathcal Z)= \pi_1[z_1]\otimes \cdots \pi_1[z_N]$ as before,
with $z_i$ formal.
Choose a partition of $N$, $\mu=(\mu_1,...,\mu_m)$, with $\mu_i$
positive integers (not necessarily ordered) and
$\mu_1+\cdots+\mu_m=N$. Let $\mathcal A = (a_1,...,a_m)$, with $a_i$
distinct complex numbers.  Let $\mathcal X = (x_1,...,x_N)$ be an
$N$-tuple of complex numbers with the first $\mu_1$ variables being
equal to $a_1$, the next $\mu_2$ equal to $a_2$, etc.

The symmetric group $S_N$ acts on $\mathcal X$ by permuting indices,
resulting in some orbit $\mathcal O_{\mathcal X}$ in $\C^N$.
The stabilizer of the $N$-tuple $\mathcal X$ is the Young subgroup
$S_\mu\in S_N$, 
$ S_\mu = S_{\{1,...,\mu_1\}}\times \cdots\times
S_{\{\mu_1+\cdots+\mu_{m-1},...,\mu_m\}}, $ where $S_{\{1,...,n\}}$ is
the group of permutations of the elements $(1,...,n)$. Corresponding
to this subgroup we have the partial Young symmetrizer in $\C[S_N]$, 
$$y_\mu = \sum_{\sigma\in S_\mu} \sigma.$$
(The tableau under
consideration is the standard one, with $(1,...,\mu_1)$ in the first
row, and so forth.)  Since $y_\mu$ is a
partial Young symmetrizer, the image of $y_\mu$ acting on
$\pi^{\otimes N}$ is, in general, not irreducible. In fact,
$$
y_\mu(\pi^{\otimes N})\simeq \pi_{\mu_1}\otimes
\pi_{\mu_2}\otimes\cdots \otimes \pi_{\mu_m}\simeq
\underset{\underset{\ell(\lambda)\leq n}{\lambda\vdash N}}{\oplus}
K_{\lambda,\mu}\ \pi_\lambda 
$$
where $K_{\lambda,\mu}$ is the Kostka number $K_{\lambda,\mu}(1)$,
and $\mu$ is the partition with parts $\mu_i$ obtained by ordering
the integers $\mu_i$.

Define the evaluation map, $\phi_{\mathcal X}$ acting on
$\C[z_1,...,z_N]$ to be the evaluation at the point $z_i=x_i$ for all
$i$:
$$
\phi_{\mathcal X} (f(z_1,...,z_N)) = f(x_1,...,x_N).
$$
Therefore the map $\nu_{\mu,\mathcal A}:\ V_N(\mathcal
Z)\to 
V_\mu(\mathcal A)=\pi_{\mu_1}[a_1]\otimes \cdots
\otimes\pi_{\mu_m}[a_m]$ is the composition of maps
\begin{equation}\label{evaluation}
{V_N(\mathcal Z)}\ {\stackrel{y_\mu\otimes\id}{\longrightarrow}} \
(\pi^{\otimes N})^{S_\mu}\otimes
\C[z_1,...,z_N]\ {\stackrel{\id\otimes\phi_{\mathcal X}}{\longrightarrow}}\
V_\mu(\mathcal A).
\end{equation}
(The order of the maps is not important as they act on different
spaces.)  Here, by $(\pi^{\otimes N})^{S_\mu}$, 
we mean the image of $y_\mu$, which
is invariant under the action of $S_\mu$.  This is a surjective map,
in particular, the second map, $\id\otimes \phi_{\mathcal X}$ is
surjective. Therefore,
\begin{equation}\label{vmu}
V_\mu(\mathcal A) \simeq (\pi^{\otimes N})^{S_\mu}\otimes \C[z_1,...,z_N]\ / \ 
{\rm ker\ }\phi_{\mathcal X}.
\end{equation}
Let $I_{\mathcal X}$ be the ideal in
$\C[z_1,...,z_N]$ of functions which vanish at the evaluation point
$z_i=x_i$. The kernel of $\phi_{\mathcal X}$ is $(\pi^{\otimes N})^{S_\mu}\otimes
I_{\mathcal X}$. Notice that $I_{\mathcal X}$ is invariant under the
action of the Young subgroup $S_\mu$.

As for the grading, one can define the grading on the space
$V_\mu(\mathcal A)$ to be the grading inherited from the preimage in
$V_N(\z)$. Therefore the isomorphism (\ref{vmu}) is an isomorphism
of filtered spaces. The associated graded space of the RHS is
isomorphic to the LHS as a graded vector space.

Next, consider the space $\F_\mu(\z)$ generated by the action of
$U(\sl_n[t])$ on the tensor product of highest weight vectors
$v_{\mu_i}$ of $\pi_{\mu_i}$. In the special case under consideration
here, with $\pi_{\mu_i}$ corresponding to the symmetric product
representation, the cyclic vector $w_\mu=v_{\mu_1}\otimes
\cdots\otimes v_{\mu_m}$ of $\F_\mu$, is the image under
$\nu_{\mu,\mathcal A}$ of $w_0 = v_0\otimes\cdots\otimes v_0\in \F_N$.
In addition, each of the modules $\pi_{\mu_i}$ is generated by the
action of the same commuting generators $b_i$ on the highest weight vector.
It follows that
$$\F_\mu(\mathcal A)= \nu_{\mu,\mathcal A}\ (\F_N(\mathcal Z)).$$
Note that the evaluation map preserves grading by definition.

\begin{lem}
Representatives of the graded tensor product $\F_\mu^*$ are a subset
of the image of representatives of $\F_{(1)^N}^*$ under the evaluation map:
$$\F_\mu^*\subset \nu_{\mu,\mathcal A}\ (\F_{(1)^N}^*)$$
where by $\F_\mu^*$ we mean the set of representatives of $\F_\mu^*$
in $F_\mu(\mathcal A)$ etc.
\end{lem}
\begin{proof}
  This follows from the definition (\ref{fusion}), because
  $\nu_{\mu,\mathcal A}(\widetilde\F_N^{(n)}) \subset
  \widetilde{\F}_\mu^{(n)}$. To see this, suppose $p^{(n)}\in
  \widetilde{\F}_N^{(n)}$. Then it can be expressed in the form
  (\ref{residues}) with $\overline{p}^{(n)}=0$.  Therefore,
$$
\nu_{\mu,\mathcal A}(p^{(n)}) = \sum_{\overset{i}{m<n}}
  f_i^{(n-m)}(x_1,...,x_N)\ \nu_{\nu,\mathcal A}(p_i^{(m)}).
  $$
  Recall that $\nu_{\nu,\mathcal A}(\F_N^{(n)}) =
  \F_\mu^{(n)}$ since the grading is preserved. Therefore $
  \nu_{\nu,\mathcal A}(p_i^{(m)})\in \F_\mu(\mathcal A)^{(m)}$ on the
  right hand side. This means that
$$\nu_{\nu,\mathcal A}(p^{(n)}) \in \C[a_1,...,a_m]\otimes
  \F_\mu(\mathcal A)^{(\leq n-1)}\cap \F_\mu(\mathcal A)^{(n)} =
  \widetilde\F_\mu(\mathcal A)^{(n)}.$$ 
  
\end{proof}

The result of Section \ref{cc} showed that $\F_{(1)^N}^*$ is the
quotient of $V_N^{S_N}$ by the ideal generated by symmetric functions
of positive degree. We will take consecutive quotients of
$V_N^{S_N}(\z)$: First by the kernel of the evaluation map and then by
$\langle J_N\rangle$.  Therefore, consider the evaluation map
$\phi_{\mathcal X}$ acting on
\begin{equation}\label{stepa}
V_N^{S_N}(\z)=
(\pi^{\otimes N} \otimes \C[z_1,...,z_N])^{S_N}.
\end{equation}

\begin{lem}
\begin{equation}
\phi_{\mathcal X}\ (\pi^{\otimes N} \otimes \C[z_1,...,z_N])^{S_N}
=\nu_{\mu,\mathcal A} \ (\pi^{\otimes N} \otimes
\C[z_1,...,z_N])^{S_N}.
\end{equation}
\end{lem}
\begin{proof}
  Consider any $w\in \pi^{\otimes N}$ and $f(z_1,...,z_N)\in
  \C[z_1,...,z_N]$.  Th set of all such $w,f$ gives a spanning set of
  vectors for $V_N(\z)^{S_N}$ after complete symmetrization:
  $$
 v = \sum_{\sigma\in S_N} \sigma(w)\otimes \sigma f(z_1,...,z_N)\in
  V_N(\z)^{S_N}.$$
Acting with $\phi_{\mathcal X}$, we have
\begin{eqnarray*}
\phi_{\mathcal X}(v) &=& \sum_{\sigma\in S_N} \sigma(w) \otimes
\phi_{\mathcal X}\sigma f(z_1,...,z_N) \\
& = & \sum_{\sigma\in S_N} \sigma(w) \otimes
\phi_{\sigma^{-1}(\mathcal X)} f(z_1,...,z_N) \\
&=& \sum_{\sigma\in S_\mu\backslash S_N}\sum_{\tau\in S_\mu}
\tau\sigma (w) \otimes \phi_{\sigma^{-1}(\mathcal X)}
f(z_1,...,z_N)
\\
&=& y_\mu \sum_{\sigma\in S_\mu\backslash S_N} \sigma(w)\otimes
\phi_{\mathcal X}\sigma f(z_1,...,z_N),
\end{eqnarray*}
(here, the sum is over coset representatives).  In the third line, we
used the fact that $(\tau\sigma)^{-1}\mathcal X = \sigma^{-1}\mathcal
X$ for $\tau\in S_\mu$. Since $y_\mu^{2}=y_\mu$, the lemma follows.
\end{proof}

Let $I_\mu$ be the symmetrization of $I_{\mathcal X}$: 
$$
I_\mu = \underset{\sigma\in S_N}{\cap} I_{\sigma(\mathcal X)}.
$$
Polynomials in $I_\mu$ are those which vanish everywhere on the
$S_N$-orbit of the point $\mathcal X$. As is clear from the proof of
the previous lemma, the kernel of $\phi_{\mathcal X}$ acting on
$V_N^{S_N}$ is the symmetrization, $(\pi^{\otimes N}\otimes I_\mu)^{S_N}$.  We
have
\begin{lem}
\begin{eqnarray}
(\pi^{\otimes N}\otimes \C[z_1,...,z_N])^{S_N}/{\rm ker}\ \phi_{\mathcal
  X} &=& 
(\pi^{\otimes N}\otimes \C[z_1,...,z_N]/I_\mu)^{S_N}\nonumber\\ 
&=& (\pi^{\otimes N}\otimes A_\mu)^{S_N}\label{quotient}
\end{eqnarray}
where $A_\mu$ is the quotient of the ring of polynomials
$\C[z_1,...,z_N]$ by the ideal of functions which vanish on the $S_N$-orbit
of $\mathcal X$.
\end{lem}

The generating polynomials of the ideal $I_\mu$ are not of homogeneous
in degree in $z_i$, and therefore the quotient (\ref{quotient}) is a
filtered space, with the filtration inherited from $V_N(\mathcal
Z)^{S_N}$.  In this picture, $a_i$ are treated as distinct complex
numbers. The graded components of (\ref{quotient}) have the same
dimensions as the graded components of the image of the evaluation
map, since it preserves grading.

  The quotient (\ref{quotient}) is isomorphic to the image of
  $V_N(\mathcal Z)^{S_N}\simeq \F_N(\mathcal Z)\otimes
  \C[z_1,...,z_N]^{S_N}$, not $\F_N(\mathcal Z)$ itself, so it is not
  equal to $\F_\mu(\mathcal A)$.  However, its associated graded space
  is the graded tensor product $\F_\mu^*$:

\begin{lem}\label{fstar}
$$\F_\mu^* \simeq \gr (\pi^{\otimes N}\otimes A_\mu)^{S_N}
$$
\end{lem}
\begin{proof}
  Taking the associated graded space is equivalent to taking the limit
  $a_i\to 0$ for all $i$. In this limit, (see Prop. 3.1 of \cite{GP})
  $I_\mu$ is identical to $J_\mu$, which contains the ideal $J_N$
  generated by symmetric functions of positive degree. Thus, taking
  the quotient by $J_N$ of (\ref{quotient}) and then taking the
  associated graded space is equivalent to simply taking the
  associated graded space.
\end{proof}

As $S_N$-modules, $A_\mu\simeq R_\mu$, and the decomposition
into irreducible $S_N$-modules has the same form as $R_\mu$:
$$
A_\mu \simeq \oplus_\lambda W_\lambda \otimes \widetilde{M}_{\lambda,\mu}.
$$
The multiplicity spaces $\widetilde{M}_{\lambda,\mu}$ are filtered
vector spaces, by degree in $z_i$. 
\begin{equation}\label{smu}
(\pi^{\otimes N}\otimes A_\mu)^{S_N} \simeq
(\pi^{\otimes N}\otimes (\underset{\lambda}{\oplus} W_\lambda \otimes
\widetilde{M}_{\lambda,\mu}))^{S_N}.
\end{equation}
As in Section \ref{cc}, we conclude that (\ref{smu}) is equal to
\begin{equation}\label{stepthree}
((\underset{\ell(\lambda)\leq n}{\oplus} \pi_\lambda \otimes W_\lambda)
\otimes (\underset{\lambda'}{\oplus} W_{\lambda'} \otimes
\widetilde{M}_{\lambda',\mu}))^{S_N} \simeq 
\underset{\overset{\lambda\vdash N}{\ell(\lambda)\leq n,\lambda\geq \mu}}
{\oplus}
\pi_\lambda\otimes \widetilde{M}_{\lambda,\mu},
\end{equation}
where the last step follows from the same reasoning as Theorem
\ref{together}. From Lemma \ref{fstar}, equation (\ref{stepthree}) and
the fact that $R_\mu=\gr A_\mu$,
\begin{cor}
\begin{equation}
  \F_\mu^* \simeq \underset{\overset{\lambda\vdash
      N}{\ell(\lambda)\leq n}}{\oplus} \pi_\lambda\otimes
  {M}_{\lambda,\mu}.
\end{equation}
\end{cor}
The restriction $\lambda\geq \mu$ turns out to be unnecessary, as the
space $M_{\lambda,\mu}$ is zero otherwise.

The main result of this section therefore is the conclusion that:
\begin{cor}
If $\mu= (\mu_1,...,\mu_m)$ is a partition of $N$ into positive
integers, then
\begin{equation}
\ch_q \F_\mu^* = \sum_{\overset{\lambda\vdash N}{\ell(\lambda)\leq
  n}}\ch \pi_\lambda
\widetilde{K}_{\lambda, \mu}(q).
\end{equation}
where $\widetilde{K}_{\lambda, \mu}(q)$ is defined in (\ref{kostilde}).
\end{cor}

\end{document}